%% file: main.tex
\documentclass[a4paper, 12pt]{amsart}
\usepackage[english]{babel}
\usepackage{amsmath, amsthm, amssymb, amsfonts, amsrefs, mathtools}
\usepackage{nccmath}
\usepackage{subfiles, hyperref}
\usepackage{xypic}
\usepackage[new]{old-arrows}
\usepackage{needspace}

\usepackage[textsize=scriptsize]{todonotes}

\include{shortcuts}

\newtheorem{definition}{Definition}
\newtheorem{lemma}[definition]{Lemma}

\newtheorem{theorem}[definition]{Theorem}
\newtheorem{corollary}[definition]{Corollary}

\theoremstyle{remark}
\newtheorem{example}[definition]{Example}
\newtheorem{remark}[definition]{Remark}

\renewcommand{\H}{\mathcal H}
\newcommand{\N}{\mathcal N}

\newcommand{\HH}{\mathbb H}

\newcommand{\GL}{\mathrm{GL}}
\newcommand{\SL}{\mathrm{SL}}

\newcommand{\Aff}{\mathrm{Aff}}
\newcommand{\SO}{\mathrm{SO}}
\newcommand{\Sym}{\mathrm{Sym}}

\newcommand{\im}{\mathrm{im}}

\newcommand{\act}[1]{\langle #1 \rangle}
\newcommand{\sla}[2]{|_{#1}\,#2}
\newcommand{\Cov}{\mathrm{Cov}}

\newcommand{\Jac}{\mathrm{Jac}}

\DeclarePairedDelimiter\floor{\lfloor}{\rfloor}

\newcommand{\rmLJ}{\rmL^\rmJ}
\newcommand{\rmRJ}{\rmR^\rmJ}
\newcommand{\wtdrmRJ}{\wtd{\rmR}\vphantom{\rmR}^\rmJ}

\newcommand{\rmRJprime}{\rmR^{\rmJ\,\prime}}

\title[Jacobi Forms of Affine Weight in Higher Cogenus]{Jacobi Forms of affine weight\\ in higher cogenus and\\ nearly holomorphic functions}
\author[J. Feldmann]{Jan Feldmann}
\author[M. Raum]{Martin Raum}
\thanks{The second named author was partially supported by Vetenskapsr\aa det Grant~2023-04217.}%

\address{%
Chalmers tekniska högskola och G\"oteborgs Universitet,
Institutionen f\"or Matematiska vetenskaper,
SE-412 96 G\"oteborg, Sweden}

\email{jan.feldmann@chalmers.se}
\email{martin@raum-brothers.eu}

\date{July, 2025}

\begin{document}

\begin{abstract}
We describe Jacobi forms of vector-valued weights in terms of classical ones, extending previous results by Ibukiyama and Kyomura to the case of arbitrary cogenus. As in their result, our isomorphisms are given by holomorphic covariant differential operators. In contrast to previous work, however, we avoid explicit calculations, which we replace by general differential geometric arguments. In the process, we obtain a structure theorem on nearly holomorphic functions on the Jacobi upper half space.
\end{abstract}

\maketitle

\subfile{sections/introduction}

\subfile{sections/praeliminaria}

\subfile{sections/nearly_holomorphic}

\subfile{sections/decomposition}

\end{document}

%% file: shortcuts.tex
\newcommand{\tx}{\text}
\newcommand{\thdash}{\nbd th}
\newcommand{\nbd}{\nobreakdash-\hspace{0pt}}

\makeatletter
\newcommand{\writelabel}[1]{#1\def\@currentlabel{#1}}
\makeatother

\newcommand{\minwidthmathbox}[2]{%
  \mathmakebox[{\ifdim#1<\width\width\else#1\fi}]{#2}%
}

\newcommand{\cN}{\mathcal{N}}

\newcommand{\rmC}{\mathrm{C}}

\newcommand{\rmJ}{\mathrm{J}}

\newcommand{\rmL}{\mathrm{L}}

\newcommand{\rmR}{\mathrm{R}}

\newcommand{\td}{\tilde}
\newcommand{\wtd}{\widetilde}
\newcommand{\ov}{\overline}

\newcommand{\wht}{\widehat}

\newcommand{\defcol}{\mathrel{:}}

\newcommand{\ra}{\rightarrow}

\newcommand{\lra}{\longrightarrow}

\newcommand{\lmto}{\longmapsto}

\newcommand{\Hom}{\operatorname{Hom}}
\newcommand{\id}{\mathrm{id}}
\renewcommand{\ker}{\operatorname{ker}}

\newenvironment{psmatrix}{\left(\begin{smallmatrix}}{\end{smallmatrix}\right)}

\newcommand{\NN}{\mathbb{N}}
\newcommand{\ZZ}{\mathbb{Z}}
\newcommand{\QQ}{\mathbb{Q}}
\newcommand{\RR}{\mathbb{R}}
\newcommand{\CC}{\mathbb{C}}

\renewcommand{\det}{\mathrm{det}}
\newcommand{\sym}{\mathrm{sym}}

\newcommand{\HS}{\mathbb{H}}

%% file: sections/introduction.tex
\section{Introduction}                                                                                 %

Jacobi forms were introduced by Eichler and Zagier~\cite{EZ}, merging the theory of elliptic functions and modular forms. Jacobi forms were pivotal in the rigorous construction of the Saito--Kurokawa lift~\cite{kurokawa-1978}, which Zagier completed~\cite{zagier-1981} following a series of papers by Maa\ss~\cite{maass-1979a,maass-1979b,maass-1979c} and Andrianov~\cite{andrianov-1979}. They continue to play an important role in the theory of Siegel modular forms based on their appearance in the  Fourier-Jacobi expansion of the latter. This connection inspired generalizations by Gritsenko and by Ziegler~\cite{Gritsenko,Ziegler}, beyond the scope of Eichler--Zagier's initial treatment, to Jacobi forms of cogenus~$h \in \NN$, which live on the Jacobi upper half space~$\HS \times \CC^h$ and have $h \times h$ matrices as their Jacobi indices, where $\HS$ denotes the complex upper half-plane.

The existing theory of Jacobi forms focuses on the scalar-valued case. However, vector-valued Siegel modular forms, such as those in Harder's Conjecture~\cite{bruinier-van-der-geer-harder-zagier-2008,atobe-chida-ibukiyama-katsurada-yamauchi-2023} or the generating series by Funke--Millson~\cite{funke-millson-2006}, motivate the study of vector-valued Jacobi forms to enable a corresponding theory of Fourier-Jacobi expansions. Such a theory would also be beneficial for employing automatic convergence statements as in~\cite{bruinier-raum-2015}. The resulting notion of vector-valued Jacobi forms, based on vector-valued weights, is distinct from the well-studied case of Jacobi forms for the Weil representation. In this paper, we offer a decomposition result for the base case of symmetric powers, thereby expressing the corresponding a priori exotic Jacobi forms in terms of classical ones.

We replace the factor of automorphy that usually appears in the transformation law of Jacobi forms~\cite{EZ,Ziegler} by one that incorporates a higher dimensional representation~$\sigma$ of~$\GL_{1+h}(\CC)$. Given a pair of elements of~$\Jac_{1, h}(\RR)$ (see~\eqref{eq:def:jacobi_group_elem_standard_notation} for notation) and~$\HS \times \CC^h$, we assign
\begin{gather*}
  \Bigl(
  \bigl( \begin{psmatrix} a & b \\ c & d \end{psmatrix}, \lambda, \mu, \kappa \bigr),\,
  (\tau, z)
  \Bigr)
\lmto
  \sigma
  \begin{psmatrix}
    c\tau + d   &   cz' - \td\lambda^\prime \\
    0           &   1
  \end{psmatrix}
\tx{.}
\end{gather*}
The right hand side typically depends on $\lambda$, unlike the classical factor of automorphy, where~$\sigma$ maps its argument to~$(c \tau + d)^k$ for some~$k \in \ZZ$, i.e., a power of the determinant.

In this paper, we focus on the case that~$\sigma$ is the~$k$\thdash{} power of the determinant times the $s$\thdash{} symmetric power representation of~$\GL_{1+h}(\CC)$ for $s \in \NN_0$, which is the most prominent case arising from our motivating connection to Siegel modular forms. Its representation space consists of degree~$s$ polynomials in~$X$ and~$Y_1, \ldots, Y_h$ and carries the linear action
\begin{gather*}
  \begin{psmatrix} r & v^\prime \\ 0 & 1 \end{psmatrix}
  f(X, Y_1, \dots, Y_h)
=
  r^k\,
  f(rX, v_1 X + Y_1, \dots, v_h X + Y_h)
\tx{,}
\end{gather*}
where~$r \in \CC^\times$, $v \in \CC^h$,~$v'$ is the transpose of~$v$, and the bottom right matrix entry~$1$ stands for the~$h \times h$ identity matrix. We give a formally complete description of the resulting notion of Jacobi forms in Definition~\ref{Jacobi_form}. The spaces of such Jacobi forms are denoted~$\rmJ_{(k,s),m}(\Gamma)$, extending the notation~$\rmJ_{k,m}(\Gamma)$ for Jacobi forms of weight~$k$ and index~$m$ for subgroups~$\Gamma \subseteq \SL_{2}(\ZZ)$. Here~$k \in \ZZ$ and~$m \in \Sym_h(\QQ)$, the space of symmetric~$h \times h$ matrices over~$\QQ$. We will assume that~$m$ is invertible and half-integral, that is, its diagonal entries lie in~$\ZZ$ and all other entries lie in~$\frac{1}{2} \ZZ$.

Connecting back to our initial discussion, we are interested in a description of~$\rmJ_{(k,s),m}(\Gamma)$ in terms of classical spaces of Jacobi forms. If~$h = 1$ such a description has been achieved by Ibukiyama and Kyomura~\cite{IK}. They proved the direct sum decomposition
\begin{gather*}
  \rmJ_{(k, s), m}(\Gamma)
\cong
  \bigoplus_{\ell = 0}^{s} \rmJ_{k + \ell, m}(\Gamma)
\tx{,}
\end{gather*}
and further established that a map from the left to the right hand side can be realized by a covariant differential operator. Note that their theorem enables a complete understanding of Fourier-Jacobi expansions of genus~$2$ Siegel modular forms in terms of classical Jacobi forms. Their work, however, relies on extensive and explicit calculations, which seem difficult to generalize to arbitrary cogenus~$h$. We employ different ideas to determine~$\rmJ_{(k,s),m}(\Gamma)$ for all~$h$, and prove the following theorem in Section~\ref{decomposition_section}.

\begin{theorem}\label{decomposition}
\Needspace*{6\baselineskip}
For $k \in \ZZ$ with $k > h \slash 2$, $s \in \NN_0$ and invertible half-integral Jacobi index~$m \in \Sym_h(\QQ)$, there is an isomorphism
\begin{gather*}
  \rmJ_{(k, s), m}(\Gamma)
\cong
\bigoplus_{\ell = 0}^{s} \rmJ_{k + \ell, m}(\Gamma)^{\oplus \binom{s - \ell + h - 1}{h - 1}}
\end{gather*}
by a covariant differential operator.
\end{theorem}

\begin{remark}
We explain the origin of the condition~$k > h \slash 2$ in Section~\ref{holomorphic_projection_subsection}.
\end{remark}

Ibukiyama--Kyomura's work is based on holomorphic differential operators, which are not necessarily covariant and thus do not preserve modularity in general. Their calculations aim at finding linear combinations of their holomorphic differential operators that do.

We reverse the approach of Ibukiyama--Kyomura, and start with differential operators that are not holomorphic but intertwine with the action of the real Jacobi group, that is, they are covariant. Such operators can be obtained by an elementary differential geometric construction (cf.\@ Section~\ref{Frobenius_reciprocity}). Their image consists of nearly holomorphic functions (cf.\@ Definition~\ref{def:nearly_holomorphic_functions}), generalizing the corresponding notion by Shimura for modular forms~\cite{Shimura1}. Our definition of nearly holomorphic functions enables and is motivated by the central Lemma~\ref{section_non_holomorphic}. We combine our covariant differential operators with a projection from nearly holomorphic functions to holomorphic ones (cf.\@ Subsection~\ref{holomorphic_projection_subsection}). The key advantage of our approach is that the construction of such holomorphic projections can be performed using general tools from differential geometry, while properties of nearly holomorphic Jacobi forms enter primarily via~Lemma~\ref{L_R_identity}. We thus avoid the analogue of the calculations performed in~\cite{IK}.

As a by-product of our work, we obtain a structure theorem for nearly holomorphic functions on $\HH\times\CC^h$ (cf.\@ Theorem~\ref{holomorphic_projection}), which is of independent interest and reminiscent of Shimura's work~\cite{Shimura1} on Hilbert modular forms and work by Pitale--Saha--Schmidt~\cite{PSS} on Siegel modular forms of genus~$2$. Notably, the proofs of neither of these related results employs differential geometric reasoning in the same way as our proofs do.

%% file: sections/praeliminaria.tex
\section{Preliminaries}\label{praeliminaria}                                                            %

We fix notation and review the standard definitions related to Jacobi forms. Throughout we let~$h \in \NN$, $k \in \ZZ$, $d, s \in \NN_0$ if not specified.

We write~$e_j$ for the~$j$\thdash{} unit vector. Denote the transpose of a matrix $m$ by $m^\prime$. Given a ring $R$, let $R^\times$ be the group of invertible
elements of $R$. Denote the set of matrices with $n$ rows and $m$ columns with coefficients in $R$ by $R^{n\times m}$,
and set
\begin{gather*}
  \Sym_n(R) = \bigl\{ m \in R^{n \times n} \mathrel{\mid} m^\prime = m \bigr\}
\tx{.}
\end{gather*}
For a square matrix $m$, set
\begin{gather*}
  e(m) = \exp\bigl( 2\pi i\, \mathrm{tr}(m) \bigr)
\end{gather*}
where $\mathrm{tr}(m)$ denotes the trace of $m$.

We write~$\Hom_\CC(V,W)$ for the space of linear maps between complex vector spaces~$V$ and~$W$.
Given a finite dimensional, complex vector space $V$ and a smooth manifold~$X$, we identify $\rmC^\infty(X)\otimes V$ with the space of smooth functions $X \to V$.

\subsection{The affine group}
Denote by $V_s$ the space of homogeneous complex polynomials of degree $s \in \NN_0$ in the indeterminates $X, Y_1, \dots, Y_h$.
We define an \emph{affine group} as
\begin{gather}
\label{eq:def:affine_group}
  \Aff_{1, h}(\CC)
=
  \bigl\{
  \begin{psmatrix} r & v^\prime \\ 0 & 1 \end{psmatrix} \in \GL_{1+h}(\CC)
  \mathrel{\mid}
  r \in \CC^\times, v \in \CC^h
  \bigr\}
\cong
  \CC^\times \ltimes \CC^h
\tx{,}
\end{gather}
where the bottom right entry~$1$ of the matrix in the middle expression is the~$h \times h$ identity matrix. This affine group acts linearly on $V_s$ by
\begin{gather*}
  \begin{psmatrix} r & v^\prime \\ 0 & 1 \end{psmatrix} f(X, Y_1, \dots, Y_h) = f(rX, v_1 X + Y_1, \dots, v_h X + Y_h)
\tx{,}
\end{gather*}
that is, we obtain a representation $(V_s, \sym^s)$ of $\Aff_{1, h}(\CC)$. In particular, $(V_0, \sym^0)$ is the trivial one-dimensional representation,
which we allow ourselves to also denote by~$\CC$. We set~$V_{-1} = \{ 0 \}$ for convenience. Another relevant representation of $\Aff_{1, h}(\CC)$ is associated with the determinant:
\begin{gather*}
  \det
\defcol
  \begin{psmatrix} r & v' \\ 0 & 1 \end{psmatrix}
\lmto
  r \in \GL_{1}(\CC)
\tx{.}
\end{gather*}
We identify the representation space of~$\det^t \otimes \sym^s$ with~$V_s$.

\begin{lemma}\label{inclusion_sym_power}
    For $s, t \in \NN_0$ with~$0 \le t \le s$, the map
    \begin{gather}
    \label{eq:def:sym_embedding}
        i_s^t \colon {\det}^{\otimes t}\otimes V_{s - t} \longrightarrow V_s, \quad f \longmapsto X^tf
    \end{gather}
    is a monomorphism of representations, and
    \begin{gather*}
        V_s \mathbin{\big\slash} \im(i_s^1) \cong \CC^{\oplus\mu(s, h)}
    \tx{,}\quad
      \tx{where }
      \mu(s, h) = \tbinom{s + h - 1}{h - 1}
    \tx{.}
    \end{gather*}
\end{lemma}

We denote the projection to the quotient in Lemma~\ref{inclusion_sym_power} and its generalization to all~$t$ by
\begin{gather}
\label{eq:def:sym_projection}
  p_s^t
\defcol
  V_s
\lra
  V_s \mathbin{\big\slash} \im(i_s^t)
\tx{.}
\end{gather}

\begin{proof}
The map~$i_s^t$ is clearly injective, and the verification
    \begin{align*}
        i_s^t\bigl(\begin{psmatrix} r & v^\prime \\ 0 & 1 \end{psmatrix} f\bigr)
            &= i_s^t\bigl(r^t\; f(rX, v_1 X + Y_1, \dots)\bigr) \\
            &= r^tX^t f(rX, v_1 X + Y_1, \dots)
             = \begin{psmatrix} r & v^\prime \\ 0 & 1 \end{psmatrix} i_s^t(f)
    \end{align*}
    is immediate. The image of $i_s^1$ consists of polynomials divisible by $X$.
    Hence, the images of the monomials in $Y_1, \dots, Y_h$ of degree $s$ under the projection
    \begin{gather}\label{projection_sym_power}
        p_s^1 \colon V_s \longrightarrow V_s \mathbin{\big\slash} \im(i_s^1)
    \end{gather}
    constitute a basis of $V_s/\im(i_s^1)$. There are $\mu(s, h)$ such monomials,
    and $\Aff_{1, h}(\CC)$ acts trivially on their image in the quotient.
\end{proof}

\subsection{Jacobi forms}

The standard references for Jacobi forms are~\cite{EZ, Ziegler}.

We refer to the positive integer~$h$ as the \emph{cogenus}, and let $R$ be a ring with identity element. The \emph{Jacobi group} $\Jac_{1, h}(R)$
consists of all matrices of the form
\begin{gather}\label{jacobi_group_elem}
  \begin{pmatrix}
  a        & 0 & b   & \td\mu'      \\
  \lambda  & 1 & \mu & \kappa       \\
  c        & 0 & d   & -\td\lambda' \\
  0        & 0 & 0   & 1            \\
  \end{pmatrix}
\in
  R^{2(1+h) \times 2(1+h)}
\end{gather}
with $\lambda, \mu \in R^h$, and
\begin{gather*}
  \begin{psmatrix} a & b \\ c & d \end{psmatrix} \in \SL_2(R)
\tx{,}\qquad
  \bigl( \td\lambda ,\, \td\mu \bigr)
  \begin{psmatrix} a & b \\ c & d \end{psmatrix}
=
  \bigl( \lambda ,\,  \mu \bigr)
\end{gather*}
and $\kappa \in R^{h\times h}$ such that $\kappa - \td\lambda\td\mu'$ is symmetric.

To match common notation, we allow ourselves to write the element in~\eqref{jacobi_group_elem} as
\begin{gather}
\label{eq:def:jacobi_group_elem_standard_notation}
  \bigl( \begin{psmatrix} a & b \\ c & d \end{psmatrix}, \lambda, \mu, \kappa \bigr)
\tx{.}
\end{gather}

The Jacobi group $\Jac_{1, h}(\RR)$ acts from the left on the \emph{Jacobi half-space} $\HH_{1, h} = \HH \times \CC^h$ by
\begin{gather*}
  g\act{\tau, z} =
  g\act{(\tau, z)} =
  \Bigl(\mfrac{a\tau + b}{c\tau + d}, \mfrac{z + \lambda \tau + \mu}{c\tau + d}\Bigr)
\end{gather*}
for $g \in \Jac_{1, h}(\RR)$ as in \eqref{jacobi_group_elem} and $(\tau, z) \in \HH_{1, h}$.

For $(\tau, z) = (x + iy, u + iv) \in \HH_{1, h}$ with $x, y\in\RR$ and $u, v \in \RR^h$, we verify that
\begin{gather}\label{transitivity}
  g = \begin{pmatrix}
    \sqrt{y}      & 0 & \frac{x}{\sqrt{y}}  & u^\prime - \frac{x}{y}v^\prime  \\
    \frac{1}{\sqrt{y}}v  & 1 & \frac{1}{\sqrt{y}}u  & 0                  \\
    0          & 0 & \frac{1}{\sqrt{y}}  & -\frac{1}{y}v^\prime        \\
    0          & 0 & 0            & 1                  \\
  \end{pmatrix}
\in
  \Jac_{1, h}(\RR)
\tx{,}
\end{gather}
and $g\act{i, 0} = (\tau, z)$. Thus $\Jac_{1, h}(\RR)$ acts transitively on $\HH_{1, h}$,
and we identify $\Jac_{1, h}(\RR)/K\cong\HH_{1, h}$ via $g \mapsto g\act{i, 0}$, where $K$ is the stabilizer of $(i, 0)$.
We identify $K$ with $\SO_2(\RR)\times \Sym_h(\RR)$. A parametrization of $K$ is given by
\begin{gather*}
  k  \colon  \RR/\ZZ \times \Sym_h(\RR)  \longrightarrow  K,      \quad
        (\theta, \kappa)      \longmapsto    k(\theta, \kappa),
\end{gather*}
where
\begin{gather}\label{k_theta_kappa}
    k(\theta, \kappa) = \begin{pmatrix}
              \cos(2\pi\theta)  & 0 & -\sin(2\pi\theta)  & 0      \\
              0          & 1 & 0          & \kappa  \\
              \sin(2\pi\theta)  & 0 & \cos(2\pi\theta)  & 0      \\
              0          & 0 & 0          & 1      \\
            \end{pmatrix}.
\end{gather}

\begin{definition}
A \emph{factor of automorphy} is a pair $(V, \eta)$ consisting of a finite-dimensional complex vector space $V$ and a smooth map
\begin{gather*}
  \eta\colon \Jac_{1, h}(\RR) \times \HH_{1, h} \to \GL(V)
\end{gather*}
that satisfies
\begin{gather*}
  \eta\bigl( g \td{g}, (\tau, z) \bigr)
=
  \eta\bigl( g, \td{g}\act{\tau, z} \bigr) \circ \eta\bigl(\td{g}, (\tau, z) \bigr)
\quad\tx{and}\quad
  \eta(\id, (\tau, z)) = \id_V
\end{gather*}
for all $g, \td{g} \in \Jac_{1, h}(\RR)$ and $(\tau, z) \in \HH_{1, h}$.
\end{definition}

Call $m \in \Sym_h(\QQ)$ \emph{half-integral} if the diagonal entries of~$m$ and all the entries of $2m$ are integers.
For $g\in \Jac_{1, h}(\RR)$ as in~\eqref{jacobi_group_elem}, $(\tau, z) \in \HH_{1, h}$, and a half-integral $m\in\Sym_h(\QQ)$, set
\begin{multline*}
  \iota_m(g, (\tau, z))
\\
\begin{aligned}
&=
  e\biggl(
    \frac{c (z + \lambda\tau + \mu)^\prime m (z + \lambda\tau + \mu)}{c\tau + d}
    - \lambda' m (2 z + \lambda \tau + \mu)
  \biggr)\,
  e(-m\kappa)
\\
&=
  e\biggl(
  \frac{(cz-\tilde{\lambda})^\prime m(\lambda\tau + z + \mu)}{c\tau + d} - \lambda^\prime m z
  \biggr)\,
  e(-m\kappa)
\end{aligned}
\end{multline*}
and
\begin{gather*}
  \eta_{(k, s), m}\bigl( g, (\tau, z) \bigr)
=
  \Bigl(
  \big({\det}^{\otimes k} \otimes \sym^s \big)
  \begin{psmatrix}
    c\tau + d   &   cz' - \td\lambda^\prime \\
    0           &   1
  \end{psmatrix}
  \Bigr)
  \otimes
  \iota_m(g, (\tau, z))
\tx{.}
\end{gather*}
A straightforward calculation confirms that $(\CC, \iota_m)$ and $(V_s, \eta_{(k, s), m})$ are factors of automorphy.
We write $\eta_{k, m}$ for $\eta_{(k, 0), m}$.

Let $(V, \eta)$ be a factor of automorphy. Given $\phi\in \rmC^\infty(\HH_{1, h})\otimes V$ and $g \in \Jac_{1, h}(\RR)$,
we obtain another function $\phi\sla{\eta}{g} \in \rmC^\infty(\HH_{1, h})\otimes V$ by
\begin{gather*}
  \phi\sla{\eta}{g} \colon
        \HH_{1, h}  \longrightarrow     V, \quad
        (\tau, z)   \longmapsto         \eta\bigl(g, (\tau, z)\bigr)^{-1}\, \phi\bigl(g\act{\tau, z}\bigr)
\tx{.}
\end{gather*}
Write $\phi\sla{(k, s), m}{g} = \phi\sla{\eta_{(k, s), m}}{g}$ and $\phi\sla{k, m}{g} = \phi\sla{(k, 0), m}{g}$.
For $g \in \Jac_{1, h}(\RR)$, the map
\begin{gather*}
    \rmC^\infty(\HH_{1, h})\otimes V \longrightarrow \rmC^\infty(\HH_{1, h})\otimes V, \quad
    \phi \longmapsto \phi\sla{\eta}{g}
\end{gather*}
defines an action of $\Jac_{1, h}(\RR)$ from the right on $\rmC^\infty(\HH_{1, h})\otimes V$.

We are now in position to define Jacobi forms with affine weight a symmetric power representation, extending the classical notion of Jacobi forms~\cite{EZ,Ziegler} and the special case~$h = 1$ by Ibukiyama--Kyomura~\cite{IK}.

\begin{definition}\label{Jacobi_form}
Let $\Gamma\subseteq \Jac_{1, h}(\ZZ)$ be a subgroup of finite index, $(V, \eta)$ a factor of automorphy and
  $m \in \Sym_h(\QQ)$ half-integral. A holomorphic function $\phi\colon \HH_{1, h} \to V$ is called a \emph{Jacobi form}
  of index $m$ with respect to $\Gamma$ and $\eta$ if:
\begin{enumerate}
\item[(i)]
for all $g \in \Gamma$, the equality $\phi\sla{\eta}{g} = \phi$ holds,

\item[(ii)]
the Fourier series expansion of~$\phi$,
\begin{gather*}
  \phi(\tau, z)
=
  \sum_{n \in \frac{1}{N} \ZZ, r \in \ZZ^h}
  \mspace{-12mu}
  c(n, r)\, e\bigl(n \tau + r^\prime z\bigr)
\tx{,}\quad
  c(n, r) \in V
\tx{,}
\end{gather*}
for suitable~$N \in \NN$ depending on~$\Gamma$, satisfy $c(n, r) = 0$ if
\begin{gather*}
  \begin{psmatrix}
    n    & r^\prime \slash 2 \\
    r \slash 2 & m
  \end{psmatrix}
\end{gather*}
is not positive semi-definite.
\end{enumerate}
\end{definition}

We write $\rmJ_{\eta, m}(\Gamma)$ for the space of all Jacobi forms of index $m$ with respect to $\Gamma$ and $\eta$.
We write $\rmJ_{(k, s), m}(\Gamma) = \rmJ_{\eta, m}(\Gamma)$ for $\eta = \eta_{(k, s), m}$.
We put $\rmJ_{k, m}(\Gamma) = \rmJ_{(k, 0), m}(\Gamma)$.

\begin{remark}\label{sum_automorphy_factors}
Let $(V, \eta)$ and $(W, \vartheta)$ be factors of automorphy. We obtain a new factor of automorphy $(V\oplus W, \eta\oplus\vartheta)$, where
\begin{gather*}
  (\eta\oplus\vartheta) \bigl(g, (\tau, z)\bigr)
=
  \eta\bigl(g, (\tau, z)\bigr) \oplus \vartheta\bigl(g, (\tau, z)\bigr)
\tx{.}
\end{gather*}
We may naturally identify
\begin{gather*}
  \rmJ_{\eta\oplus\vartheta, m}(\Gamma)
\cong
  \rmJ_{\eta, m}(\Gamma) \oplus \rmJ_{\vartheta, m}(\Gamma)
\tx{.}
\end{gather*}
\end{remark}

%% file: sections/nearly_holomorphic.tex
\section{Nearly holomorphic functions}\label{nearly_holomorphic}                                        %

A standard reference for nearly holomorphic functions, also called almost holomorphic functions, is~\cite{Shimura1}.

Throughout this section, we assume that~$\nu \in \NN_0^h$, $r \in \NN_0$, and~$m \in \Sym_h(\QQ)$ half-integral, if not specified.

In this section, we introduce nearly holomorphic functions on $\HH_{1, h}$ and prove
Theorem~\ref{holomorphic_projection}, which provides a structural result similar to
\cite[Theorem 5.2]{Shimura1}, though less conclusive.
For $(\tau, z) = (x + iy, u + iv) \in \HH_{1, h}$ with real $x, y \in \RR$, $u, v \in \RR^h$, set
\begin{gather*}
    \alpha  =   \left(\alpha_1,          \ldots, \alpha_h         \right)^\prime
            =   \left(\mfrac{v_1}{y},    \ldots, \mfrac{v_h}{y}   \right)^\prime,
    \qquad
    \beta   =   \mfrac{1}{8\pi y},
\end{gather*}
and put
\begin{gather*}
    \alpha^\nu = \alpha_1^{\nu_1} \cdots \alpha_h^{\nu_h}.
\end{gather*}

Write $\H_{(k, s), m}$ for the space of holomorphic functions $\HH_{1, h} \to V_s$, which the Jacobi group acts on by the slash action~$|_{(k,s), m}$.
We let $\H_{(k, s), m}[\alpha, \beta]$ denote the ring of polynomials in $\alpha_1, \ldots, \alpha_h$ and $\beta$ with
coefficients in~$\H_{(k, s), m}$. Elements of $\H_{(k, s), m}[\alpha, \beta]$ are called \emph{nearly holomorphic} functions
on $\HH_{1, h}$ with values in~$V_s$. By choosing a basis $b_1, \dots, b_n$ of $V_s$, any such nearly holomorphic
$f \colon \HH_{1, h} \to V_s$ may be written as
\begin{gather}\label{nearly_holomorphic_f}
  f
=
  \sum_{\nu, r} \alpha^\nu \beta^r\, f_{\nu,r}
\quad\tx{with }
  f_{\nu,r}
=
  \sum_{i = 1}^n
  f_{\nu, r}^i\, b_i
\tx{.}
\end{gather}

Set
\begin{gather*}
  |\nu, r| = \nu_1 + \cdots + \nu_h + 2r
\tx{,}
\end{gather*}
and define the \emph{degree} of $f \neq 0$ as in \eqref{nearly_holomorphic_f} by
\begin{gather*}
  \deg f = \max \bigl\{ |\nu, r| \mathrel{\mid} f_{\nu, r} \neq 0 \bigr\}
\tx{,}\quad\tx{and }
  \deg 0 = -\infty
\tx{.}
\end{gather*}
Also for~$f \ne 0$ set
\begin{gather*}
  \deg(\alpha_j, f) = \max \bigl\{ \nu_j \mathrel{\mid} f_{\nu, r} \neq 0 \bigr\}
\tx{,}\quad
  \deg(\beta, f) = \max \bigl\{ r \mathrel{\mid} f_{\nu, r} \neq 0 \bigr\}
\tx{,}
\end{gather*}
and define
\begin{gather*}
  \deg(\alpha_j, 0) = \deg(\beta, 0) = -\infty
\tx{.}
\end{gather*}

\begin{definition}
\label{def:nearly_holomorphic_functions}
The space of depth~$d$ nearly holomorphic functions for weight~$(k,s)$ is
\begin{gather*}
  \N^d_{(k, s), m}
=
  \bigl\{
  f \in \H_{(k, s), m}[\alpha, \beta]
  \mathrel{\mid}
  \deg\bigl( p^t_s\circ f \bigr) < d + t\ \tx{ for all } t \leq s
  \bigr\}
\tx{,}
\end{gather*}
where the projections~$p^t_s$ are given in~\eqref{eq:def:sym_projection}.
\end{definition}

Note that~$\cN^d_{(k,s),m}$ is independent of~$k$ and~$m$ as a set,
but the action of the Jacobi group on it depends on them.

\begin{example}
The function~$\HS_{1, 1} \ra V_2$ given by $(\alpha^2 + \beta)X^2 + \alpha XY + Y^2$ has depth~$0$.
\end{example}

For brevity, write
\begin{gather*}
    \H_{k, m} = \H_{(k, 0), m}, \qquad \N^d_{k, m} = \N^d_{(k, 0), m}.
\end{gather*}

\subsection{Differential operators}
We recall the Maa{\ss}-Shimura raising and lowering operators for Jacobi forms.
We refer to~\cite{CR} for information on how to construct these operators. We put
\begin{gather}
\label{eq:def:renomalized_diff}
\begin{alignedat}{3}
  \partial_\tau   &= \mfrac{1}{4\pi i} \mfrac{\partial}{\partial\tau}
\tx{,}\quad
&
  \partial_{z, j} &= \mfrac{1}{4\pi i} \mfrac{\partial}{\partial z_j}
\tx{,}\quad
&
  \partial_z      &= \bigl(\partial_{z, 1},\dots, \partial_{z, h}\bigr)^\prime
\tx{,}
\\
  \partial_{\ov{\tau}} &= \mfrac{1}{4\pi i} \mfrac{\partial}{\partial\overline{\tau}}
\tx{,}\quad
&
  \partial_{\ov{z}, j} &= \mfrac{1}{4\pi i} \mfrac{\partial}{\partial\overline{z}_j}
\tx{,}\quad
&
  \partial_{\ov{z}}    &= \bigl(\partial_{\ov{z}, 1},\dots, \partial_{\ov{z}, h}\bigr)^\prime
\tx{,}
\end{alignedat}
\end{gather}
and
\begin{gather}
\label{eq:def:renomalized_var}
  \hat{y}   = 8 \pi y
\tx{,}\quad
  \hat{v}_j = 8 \pi v_j
\tx{,}\quad
  \hat{v}   = \left(\hat{v}_1, \dots, \hat{v}_h\right)^\prime
\tx{.}
\end{gather}
We consider the following differential operators:
\begin{gather}
\label{eq:def:maass_operators}
\begin{aligned}
  \rmR_k
&=
  \partial_\tau + \alpha^\prime \partial_z + \tfrac{1}{2} \alpha^\prime m \alpha - k\beta
\tx{,}
&
  \rmRJ
&=
  \partial_z + m\alpha
\tx{,}
\\
  \rmL
&=
  \hat{y}\left(\hat{y}\partial_{\ov{\tau}} + \hat{v}^\prime\partial_{\ov{z}}\right)
\tx{,}
&
  \rmLJ
&=
  \hat{y} \partial_{\ov{z}}
\tx{.}
\end{aligned}
\end{gather}
Note that we suppress the dependence of the raising operators on the Jacobi index~$m$ from our notation. We write~$\rmRJ_j$ and~$\rmLJ_j$ for the~$j$\thdash{} component of~$\rmRJ$ and~$\rmLJ$.

\begin{remark}
\label{rm:maass_operators_previous_definition}
We renormalize the operators in \cite[Definition 2.5]{CR}:
\begin{gather*}
  \rmR_k = -\mfrac{1}{8\pi}\, X^{k, L}_+
\tx{,}\quad
  \rmRJ  = -\mfrac{1}{4\pi}\, Y^{k, L}_+
\tx{,}\quad
  \rmL   = 8\pi\, X_-^{k, L}
\tx{,}\quad
  \rmLJ  = 2\, Y_-^{k, L}
\tx{.}
\end{gather*}
In the above equations our operators appear on the left-hand side while the right-hand side is expressed
in terms of the operators found in the reference.
\end{remark}

\begin{remark}\label{derivative_alpha_beta}
The operators $\rmL$ and $\rmLJ_j$ act on nearly holomorphic functions similar to the derivative
in $\alpha_j$ and $\beta$ respectively, that is,
\begin{gather*}
  \rmL\bigl( \alpha^\nu \beta^r f \bigr)
=
  -r\, \alpha^\nu \beta^{r-1}\, f
\tx{,}\quad
  \rmLJ_j\bigl(\alpha^\nu \beta^r f \bigr)
=
  \nu_j\, \alpha^{\nu - e_j^\prime} \beta^r\, f
\end{gather*}
for a holomorphic function $f \colon \HH_{1, h} \to \CC$.
\end{remark}

\begin{remark}\label{commutator_1}
We have~$[\rmLJ_i, \rmLJ_j] = [\rmRJ_i, \rmRJ_j ] = 0$ for all~$i,j$ by~\cite[Proposition 2.6]{CR}, and the non-zero commutator relations
\begin{align*}
  \left[\rmL,    \rmR_k \right] = k
\tx{,}\quad
  \left[\rmL,    \rmRJ_j\right] = \rmLJ_j
\tx{,}\quad
  \left[\rmLJ_j, \rmR_k \right] = \rmRJ_j
\tx{,}\quad
  \left[\rmLJ_i, \rmRJ_j\right] = m_{i, j}
\tx{,}
\end{align*}
where we adjust the weight for the raising operator according to its covariance properties as follows:
\begin{gather*}
  \left[\rmL, \rmR_k\right] = \rmL \circ \rmR_k - \rmR_{k-2} \circ \rmL
\quad\tx{and}\quad
  \left[\rmLJ, \rmR_k\right]= \rmLJ \circ \rmR_k - \rmR_{k-1} \circ \rmLJ
\tx{.}
\end{gather*}
\end{remark}

We are interested in constructing holomorphic covariant operators, which we can then restrict to spaces of Jacobi forms. As an intermediate step, we consider
\begin{gather}
\label{eq:def:differential_operators_vy_coefficients}
  \CC\big[ \alpha, \beta \big]
\otimes
  \CC\big[ \partial_\tau, \partial_{\ov\tau}, \partial_z, \partial_{\ov{z}} \big]
\tx{,}
\end{gather}
the space of differential operators on~$\HS_{1,h}$ whose coefficients are polynomials in~$\alpha$ and~$\beta$. It is a (non-commutative) algebra by
\begin{gather}
\label{eq:renomalized_diff_action}
  \partial_\tau\, \alpha_j       = \alpha_j \beta
\tx{,}\quad
  \partial_{z,j}\, \alpha_j      = - \beta
\tx{,}\quad
  \partial_{\ov\tau}\, \alpha_j  = - \alpha_j \beta
\tx{,}\quad
  \partial_{\ov{z},j}\, \alpha_j = \beta
\tx{.}
\end{gather}
We define the order of an element of~\eqref{eq:def:differential_operators_vy_coefficients} as its degree in the second tensor component, that is, in~$\partial_\tau, \partial_{\ov\tau}, \partial_z, \partial_{\ov{z}}$.

\begin{definition}
\label{covariant_diffop}
Let $(V, \eta)$ and $(W, \vartheta)$ be two factors of automorphy. We call a map
\begin{gather*}
  D
\defcol
  \rmC^\infty (\HH_{1, h}) \otimes V \longrightarrow \rmC^\infty(\HH_{1, h}) \otimes W
\end{gather*}
a \emph{covariant differential operator} if
\begin{gather}
\label{eq:def:covariant_diffop:coefficients}
  D
\in
  \Big(
  \CC\big[ \alpha, \beta \big]
  \otimes
  \CC\big[ \partial_\tau, \partial_{\ov\tau}, \partial_z, \partial_{\ov{z}} \big]
  \Big)
  \,\otimes\,
  \Hom_\CC(V,W)
\end{gather}
and for all $f \in \rmC^{\infty}(\HH_{1, h}) \otimes V$ and all~$g \in \Jac_{1, h}(\RR)$ we have
\begin{gather}
\label{eq:def:covariant_diffop:covarianc}
  D\big( f\sla{\eta}{g} \big)
=
  D(f)\sla{\vartheta}{g}
\tx{.}
\end{gather}
\end{definition}

We let~$\Cov(\eta, \vartheta)$ denote the space of covariant differential operator from~$\eta$ to~$\vartheta$.

Given~$k' \in \ZZ$ and~$d', s' \in \NN_0$, we consider the factors of automorphy~$\eta = \eta_{(k,s),m}$ and~$\vartheta = \eta_{(k',s'), m}$. We say that~$D \in \Cov(\eta, \vartheta)$ is a covariant differential operator for nearly holomorphic or for holomorphic functions if it restricts to a map
\begin{gather*}
  \N^d_{(k, s), m} \lra \N^{d'}_{(k', s'), m}
\quad\tx{or}\quad
  \H_{(k, s), m}   \lra \H_{(k', s'), m}
\tx{.}
\end{gather*}

The operators in~\eqref{eq:def:maass_operators} are covariant differential operators of order one: They satisfy~\eqref{eq:def:covariant_diffop:coefficients} by inspection of the defining expressions and are covariant as in~\eqref{eq:def:covariant_diffop:covarianc} by~\cite[Proposition 2.6]{CR} in light of Remark~\ref{rm:maass_operators_previous_definition}. We next showcase the simplest differential operators of order zero.

\begin{example}
The inclusion $i^t_s \colon {\det}^{\otimes t} \otimes V_{s-t} \to V_s$ in~\eqref{eq:def:sym_embedding} yields by composition a covariant differential operator of order zero
\begin{gather}
\label{eq:def:sym_embedding_diffop}
  \big( i^t_s \big)_\ast
\defcol
  \N^{d+t}_{(k+t, s-t), m}
\lra
  \N^{d}_{(k, s), m}
\tx{,}\quad
  f
\lmto
  i^t_s \circ f
\tx{.}
\end{gather}
Moreover, the short exact sequence
\begin{gather*}
    0 \longrightarrow {\det} \otimes V_{s-1} \longrightarrow V_s \longrightarrow V_s \mathbin{\slash} {\im(i^1_s)} \cong \CC^{\oplus\mu(s, h)} \longrightarrow 0
\end{gather*}
with inclusion~$i^1_s$ and projection~$p^1_s$ from \eqref{eq:def:sym_embedding} and~\eqref{eq:def:sym_projection} and~$\mu(s,h)$
as in Lemma~\ref{inclusion_sym_power} yields by composition the two short exact sequences whose morphisms are covariant differential operators:
\begin{gather}\label{holomorphic_short_exact}
    0 \longrightarrow \H_{(k+1, s-1), m} \longrightarrow \H_{(k, s), m} \longrightarrow \H_{k, m}^{\oplus \mu(s, h)} \longrightarrow 0
\end{gather}
and
\begin{gather}\label{nearly_holomorphic_short_exact}
    0 \longrightarrow \N^{d+1}_{(k+1, s-1), m} \longrightarrow \N^{d}_{(k, s), m} \longrightarrow \bigl(\N^{d}_{k, m}\bigr)^{\oplus \mu(s, h)} \longrightarrow 0
\tx{.}
\end{gather}
\end{example} 

\begin{lemma}\label{diffop_and_fourier_expansion}
For~$k' \in \ZZ$ and~$s' \in \NN_0$, any covariant operator for holomorphic functions
    \begin{gather*}
        D \colon \H_{(k, s), m} \longrightarrow \H_{(k', s'), m}
    \end{gather*}
    restricts to a map
    \begin{gather*}
        \rmJ_{(k, s), m}(\Gamma) \longrightarrow \rmJ_{(k', s'), m}(\Gamma).
    \end{gather*}
\end{lemma}
\begin{proof}
    For $g \in \Gamma$ and $\phi \in \rmJ_{(k, s), m}(\Gamma)$, the calculation
    \begin{gather*}
        D(\phi)\sla{(k', s'), m}{g} = D(\phi\sla{(k, s), m}{g}) = D(\phi)
    \end{gather*}
    shows that $D(\phi)$ satisfies the first condition of Definition~\ref{Jacobi_form}.
    Furthermore, since the Fourier expansion of~$\phi$ converges absolutely and locally uniformly,
    we may differentiate it term by term to obtain a Fourier expansion of~$D(\phi)$.
    Thus the second condition of Definition~\ref{Jacobi_form} is satisfied as well and $D(\phi) \in \rmJ_{(k', s'), m}(\Gamma)$.
\end{proof}

\begin{corollary}%
\label{cor:jacobi_short_left_exact_sequence}
The short exact sequence~\eqref{holomorphic_short_exact} restricts to an exact sequence
\begin{gather}
\label{eq:cor:jacobi_short_left_exact_sequence}
    0 \longrightarrow \rmJ_{(k+1, s-1), m}(\Gamma) \longrightarrow \rmJ_{(k, s), m}(\Gamma) \longrightarrow \rmJ_{k, m}(\Gamma)^{\oplus \mu(s, h)}.
\end{gather}
\end{corollary}

\begin{remark}\label{non-covariant_split}
    The pushforward along any $\CC$\nbd{}linear section to the morphism~$p^1_s\colon V_s \to V_s \mathbin{\slash} \im(i^1_s)$
    yields a splitting of the short exact sequences~\eqref{holomorphic_short_exact} and~\eqref{nearly_holomorphic_short_exact}.
    However, such a splitting is not covariant in general and therefore cannot be restricted to the corresponding spaces of Jacobi forms.
    
    A priori, it is hence not clear whether the final morphism in~\eqref{eq:cor:jacobi_short_left_exact_sequence} is surjective.
\end{remark}

\subsection{A holomorphic projection}\label{holomorphic_projection_subsection}
We now investigate the interaction of the covariant differential operators introduced above with nearly holomorphic functions and prove our main result on nearly holomorphic functions on $\HH_{1, h}$, Theorem~\ref{holomorphic_projection}. Lemma~\ref{L_R_identity} contains a key computation of eigenvalues of covariant differential operators on holomorphic functions, which gives rise to the assumption~$k > h \slash 2$ in both Theorems~\ref{decomposition} and~\ref{holomorphic_projection}.

Assume that~$m$ is invertible. It is convenient to separate contributions of the monomials $\alpha_1, \ldots, \alpha_h$ and $\beta$ to the raising operators. Set
\begin{alignat*}{2}
  \wtd\Delta_k &= \rmR_k - \tfrac{1}{2} \rmRJprime m^{-1} \rmRJ
              &&= \partial_\tau - \tfrac{1}{2} \partial_z^\prime m^{-1} \partial_z + \left(\tfrac{h}{2} - k\right) \beta
\tx{,}
\\
  \wtdrmRJ &= m^{-1}\, \rmRJ
\tx{,}
\quad\text{i.e., }
  \wtdrmRJ_{j} &&= \sum_{i = 1}^h (m^{-1})_{j, i}\, \partial_{z, i} + \alpha_j
\tx{.}
\end{alignat*}
Note that the first operator is an analogue of the heat operator.

\begin{remark}\label{commutator_2}
The following commutator relations for~$i \ne j$ can be deduced from Remark~\ref{commutator_1}:
\begin{align*}
  \bigl[\rmL, \wtd\Delta_k\bigr] &= k - \tfrac{h}{2} - \tfrac{1}{2} \rmRJprime m^{-1} \rmLJ
\tx{,}
\\[.1\baselineskip]
  \bigl[\rmRJ_j, \wtdrmRJ_j\bigr] &=
  \bigl[\rmRJ_i, \wtdrmRJ_j\bigr] = 0
\tx{,}\quad
  \bigl[\rmLJ_j, \wtdrmRJ_j\bigr] = 1
\tx{,}\quad
  \bigl[\rmLJ_i, \wtdrmRJ_j\bigr] = 0
\tx{,}
\end{align*}
where we adjust the weight of~$\wtd\Delta_k$ in the first commutator as in Remark~\ref{commutator_1}.
\end{remark}

In the next two lemmas, we record the interaction between degrees and Maa\ss{}--Shimura operators in cases that we will employ later.

\begin{lemma}\label{degree_raising}
    Let $f \colon \HH_{1, h} \to \CC$ be a nearly holomorphic function. Then for all~$j$ and all~$i \ne j$
    \begin{align*}
        \deg(\alpha_j, \wtdrmRJ_j(f)) &= \deg(\alpha_j, f) + 1\tx{,} \\
        \deg(\alpha_i, \wtdrmRJ_j(f)) &\leq \deg(\alpha_i, f)\tx{,}\quad
        \deg(\beta, \wtdrmRJ_j(f)) = \deg(\beta, f)
      \tx{.}
    \end{align*}
    Further, for all~$j$
    \begin{gather*}
      \deg(\alpha_j, \wtd\Delta_k(f)) \leq \deg(\alpha_j, f)
    \tx{,} 
    \end{gather*}
    and if $\deg(\alpha_j, f) = 0$ for all $j$, then
    \begin{gather*}
      \deg(\beta, \wtd\Delta_k(f)) = \deg(\beta, f) + 1
    \tx{.}
    \end{gather*}
\end{lemma}
\begin{proof}
    For a holomorphic function
    $f \colon \HH_{1, h} \to \CC$, calculate
    \begin{gather*}
        \wtdrmRJ_j\big( \alpha^\nu \beta^r f \big)
     =
       \sum_{i = 1}^h (m^{-1})_{j, i}\,
       \Bigl(\nu_i \alpha^{\nu - e_i^\prime} \beta^{r + 1} f
               + \alpha^\nu\beta^r \partial_{z, i} f\Bigr)
           + \alpha^{\nu + e_j^\prime} \beta^r f,
    \end{gather*}
    which demonstrates the first two claims. The remaining calculations are similar and therefore omitted.
\end{proof}

\begin{lemma}\label{degree_lowering}
    Consider a nearly holomorphic function $f \colon \HH_{1, h} \to \CC$. If\/~$\deg(\alpha_j, f) \leq 0$, we have~$\rmLJ_j(f) = 0$, and otherwise for all~$i \neq j$
    \begin{align*}
      \deg(\alpha_j, \rmLJ_j(f)) &= \deg(\alpha_j, f) - 1\tx{,}\\
      \deg(\alpha_i, \rmLJ_j(f)) &\le \deg(\alpha_i, f)\tx{,}\quad
      \deg(\beta, \rmLJ_j(f)) \le \deg(\beta, f)
    \tx{.}
    \end{align*}
    If\/~$\deg(\beta, f) \leq 0$, we have~$\rmL(f) = 0$, and otherwise for all~$j$
    \begin{align*}
      \deg(\alpha_j, \rmL(f)) &\le \deg(\alpha_j, f)\tx{,}\\
      \deg(\beta, \rmL(f)) &= \deg(\beta, f) - 1
    \tx{.}
    \end{align*}
\end{lemma}
\begin{proof}
    The lemma follows from the formulas in Remark~\ref{derivative_alpha_beta}.
\end{proof}

\begin{corollary}\label{covariant_operators}
The operators $\rmLJ_j$ and $\rmL$ restrict to covariant operators
\begin{align*}
  \wtd\Delta_{k} \colon  \N^{d}_{k, m}   &\lra \N^{d+2}_{k+2, m},     &
  \wtdrmRJ_j     \colon  \N^{d}_{k, m}   &\lra \N^{d+1}_{k+1, m},     \\
  \rmL           \colon  \N^{d}_{k, m}   &\lra \N^{d-2}_{k-2, m},     &
  \rmLJ_j        \colon  \N^{d}_{k, m}   &\lra \N^{d-1}_{k-1, m}.
\end{align*}
\end{corollary}
Note that by our definition, nearly holomorphic functions of negative depth are zero.
\begin{proof}
For the covariance of the operators, we refer to Remark~\ref{rm:maass_operators_previous_definition} and~\cite[Proposition~2.6]{CR}. All required inequalities for the degree are given in Lemmas~\ref{degree_raising} and~\ref{degree_lowering}.
\end{proof}

Define the following operators on nearly holomorphic functions:
\begin{gather}
\label{eq:def:L_R_for_section}
\begin{aligned}
    \wht\rmR_{\nu, r}   &=  (\wtdrmRJ_1)^{\nu_1} \circ \cdots \circ (\wtdrmRJ_h)^{\nu_h}
                        \circ \wtd\Delta_{k - |\nu, r| + 2(r-1)} \circ \cdots \circ \wtd\Delta_{k - |\nu, r|}
  \tx{,}
  \\
    \wht\rmL_{\nu, r}   &=  \rmL^r \circ (\rmLJ_h)^{\nu_h} \circ \cdots \circ (\rmLJ_1)^{\nu_1}
  \tx{,}
\end{aligned}
\end{gather}
where we omit the dependence on~$k$ and~$m$.
By Remark~\ref{commutator_1} and~\ref{commutator_2}, we may reorder the operators $\wtdrmRJ_1, \ldots, \wtdrmRJ_h$
and $\rmLJ_1, \ldots, \rmLJ_h$ ad libitum in the composition above. This fact will be used tacitly in
the calculations below.

By Corollary~\ref{covariant_operators}, the above operators restrict to covariant operators
\begin{gather*}
    \wht\rmR_{\nu, r}   \colon \H_{k-|\nu, r|, m} \longrightarrow \N^{|\nu, r|}_{k, m}
\quad\tx{and}\quad
    \wht\rmL_{\nu, r}   \colon \N^{|\nu, r|}_{k, m} \longrightarrow \H_{k-|\nu, r|, m}
\tx{.}
\end{gather*}
More specifically, by Lemma~\ref{degree_raising} if $f \colon \HH_{1, h} \to \CC$ is holomorphic, we have the inequality
\begin{gather}
\label{degree_R}
    \deg(\alpha_j, \wht\rmR_{\nu, r}(f)) \leq \nu_j
\end{gather}
for all $j$.

\begin{lemma}\label{L_R_identity}
    If $k - d > h \slash 2$ for~$d = |\nu, r|$, then the composite
    \begin{gather*}
        \wht\rmL_{\nu, r} \circ \wht\rmR_{\nu, r} \colon \H_{k-d, m} \longrightarrow \H_{k-d, m}
    \end{gather*}
    is a positive multiple of the identity.
\end{lemma}

\begin{proof}
We proceed by induction on $d$, where the case~$d = 0$ is vacuous.

First, we consider the case~$\nu = 0$. We let~$f \colon \HS_{1,h} \ra \CC$ be holomorphic. Combining Lemma~\ref{degree_raising} and~\ref{degree_lowering}, we find that for all~$1 \le q \le r$
\begin{gather*}
  \rmLJ \circ \wtd\Delta_{k-2q} \circ \cdots \circ \wtd\Delta_{k-2r} (f)
=
  0
\tx{.}
\end{gather*}
The commutator relations in Remark~\ref{commutator_2}, which a priori give rise to a contribution of~$\rmLJ$, thus yield
\begin{align*}
  \wht\rmL_{0, r} \circ \wht\rmR_{0, r}(f)
&=
  \rmL^r \circ \wtd\Delta_{k-2} \circ \cdots \circ \wtd\Delta_{k-2r} (f)
\\
&=
  \Bigl( \sum_{q = 1}^r \bigl( k-2q - \tfrac{h}{2} \bigr)\Bigr)\,
  \rmL^{r-1} \circ \wtd\Delta_{k-4} \circ \cdots \circ \wtd\Delta_{k-2r} (f)
\tx{.}
\end{align*}
Since~$d = 2r$, the assumption~$k - d > h \slash 2$ ensures that the coefficient on the right hand side is positive, and the lemma follows by induction on~$d$.

Next, we consider the case~$\nu \neq 0$. Let $j$ be the smallest index such that $\nu_j \neq 0$.
We record that if~$f \in \HH_{1, h} \to \CC$ is nearly holomorphic with $\deg(\alpha_j, f) = 0$, then~$\rmLJ_j(f) = 0$ and thus
\begin{gather}
\label{eq:prf:la:L_R_identity:auxiliary}
  (\rmLJ_j)^{\nu_j} \circ (\wtdrmRJ_j)^{\nu_j}\, (f)
=
  \nu_j\, (\rmLJ_j)^{\nu_j - 1}\circ (\wtdrmRJ_j)^{\nu_j - 1}\, (f)
\end{gather}
by repeated application of the commutator relations in Remark~\ref{commutator_2}.

Now let $f \colon \HH_{1, h} \to \CC$ be holomorphic. We can apply the relation in~\eqref{eq:prf:la:L_R_identity:auxiliary} to~$\wht\rmR_{\nu - \nu_j e_j^\prime, r}(f)$, since its degree in~$\alpha_j$ is~$0$ by~\eqref{degree_R}. We conclude that
\begin{align*}
  \wht\rmL_{\nu, r} \circ \wht\rmR_{\nu, r}\, (f)
&=
  \wht\rmL_{\nu-\nu_j e_j^\prime,r}\circ (\rmLJ_j)^{\nu_j}\circ (\wtdrmRJ_j)^{\nu_j}\circ \wht\rmR_{\nu-\nu_j e_j^\prime,r}\, (f) \\
&=
  \nu_j\, \wht\rmL_{\nu - e_j^\prime, r}\circ \wht\rmR_{\nu - e_j^\prime, r}\, (f)
\tx{,}
\end{align*}
which is a positive multiple of~$f$ by the induction hypothesis.
\end{proof}

\begin{lemma}\label{L_R_image_kernel}
Given distinct pairs $(\mu, q)$ and $(\nu, r)$ in~$\NN_0^h \times \NN_0$ with $|\mu, q| = |\nu, r| = d$ and $q \leq r$, for all holomorphic functions~$f$ on~$\HS_{1,h}$, we have
\begin{gather*}
  \wht\rmL_{\mu, q} \circ \wht\rmR_{\nu, r}\, (f)
=
  0
\tx{.}
\end{gather*}
\end{lemma}
\begin{proof}
    Note that $\nu_j < \mu_j$ for some index $j$. Thus
    \begin{gather*}
        \deg\bigl(\alpha_j, \wht\rmR_{\nu, r}(f)\bigr) < \mu_j
    \end{gather*}
    for holomorphic $f$ by~\eqref{degree_R}. Using Lemma~\ref{degree_lowering}, we calculate
    \begin{gather*}
        \wht\rmL_{\mu, q} \circ \wht\rmR_{\nu, r} (f) = \wht\rmL_{\mu - \mu_j e_j^\prime, q} \circ (\rmLJ_j)^{\mu_j} \circ \wht\rmR_{\nu, r}(f) = 0,
    \end{gather*}
    which proves the claim.
\end{proof}

\begin{theorem}\label{holomorphic_projection}
If $k - d > h \slash 2$ and~$m$ is invertible, we have the direct sum decomposition
\begin{gather*}
    \N^d_{k, m} \cong \bigoplus_{\ell = 0}^d \H_{k - \ell, m}^{\oplus n(\ell)}
\end{gather*}
by a covariant differential operator, where
\begin{gather*}
  n(\ell) = \# \bigl\{ (\nu, r) \in \NN_0^h \times \NN_0 \mathrel{\mid} |\nu, r| = \ell \bigr\}
\tx{.}
\end{gather*}
\end{theorem}
\begin{proof}
Denote by $D_{-1}$ the zero map $\N^d_{k, m} \to 0$. For $q \in \NN_0$ with $0 \leq q \leq d/2$, define inductively the operator
\begin{gather*}
    D_q \colon \ker(D_{q-1}) \longrightarrow \bigoplus_{|\nu, q| = d} \H_{k - d, m},
    \quad f \longmapsto \bigoplus_{|\nu, q| = d} \wht\rmL_{\nu, q}(f),
\end{gather*}
where the direct sum runs over all $\nu \in \NN_0^h$ with $|\nu, q| = d$. This map is well-defined by a calculation invoking Lemma~\ref{degree_lowering}.

Since $k - d > h \slash 2$, Lemma~\ref{L_R_identity} applies and we find scalars $c_{\nu, r}$ such that $c_{\nu, r} \wht\rmR_{\nu, r}$
is a right inverse to $\wht\rmL_{\nu, r}$ on holomorphic functions.
Combining this with Lemma~\ref{L_R_image_kernel}, we deduce that we have an operator
\begin{gather*}
  S_q \colon 
  \bigoplus_{|\nu, q| = d} \H_{k - d, m} \longrightarrow \ker(D_{q-1}),
\quad
  (f_{\nu,q})_{(\nu,q)} \longmapsto \mspace{-2mu}\sum_{|\nu, q| = d}\mspace{-4mu} c_{\nu, q}\, \wht\rmR_{\nu, q}(f_{\nu,q})
\tx{,}
\end{gather*}
which is a section to~$D_q$.
Thus we obtain a split exact sequence
\begin{gather*}
    0 \longrightarrow \ker(D_{q})
    \longrightarrow \ker(D_{q-1})
    \longrightarrow \bigoplus_{|\nu, q| = d} \H_{k - d, m} \longrightarrow 0,
\end{gather*}
in which all maps are covariant differential operators.

Another calculation with monomials in~$\alpha_j$ and~$\beta$ shows that
\begin{gather*}
    \ker\bigl(D_{\floor{\frac{d}{2}}}\bigr) = \bigcap_{|\nu, r| = d} \ker\bigl(\wht\rmL_{\nu, r}\bigr) = \N^{d-1}_{k, m}.
\end{gather*}
We conclude that
\begin{gather*}
    \N^d_{k, m} \cong \N^{d-1}_{k, m} \oplus \bigoplus_{|\nu, r| = d} \H_{k - d, m}
\end{gather*}
by a covariant operator, and the theorem follows by induction on~$d$.
\end{proof}

%% file: sections/decomposition.tex
\section{Decomposition of vector-valued Jacobi forms}                                                   %
In this section, we prove our main result, Theorem~\ref{decomposition}, which decomposes vector-valued Jacobi forms
and extends \cite[Theorem~1.1]{IK} to higher cogenus.

We maintain the assumption~$m \in \Sym_h(\QQ)$ half-integral from Section~\ref{nearly_holomorphic}.

We return to the short exact sequence
\begin{gather*}
    0   \longrightarrow             \rmJ_{(k + 1, s - 1), m}(\Gamma)
        \longrightarrow             \rmJ_{(k, s), m}(\Gamma)
        \longrightarrow             \rmJ_{k, m}(\Gamma)^{\oplus \mu(s, h)}
\end{gather*}
in Corollary~\ref{cor:jacobi_short_left_exact_sequence}. Our goal is to show that the right arrow is surjective, and we will achieve this by showing that it splits. In Section~\ref{Frobenius_reciprocity}, we split the analogous short exact sequence of spaces of nearly holomorphic functions.
Then, in Section~\ref{decomposition_section}, we use the results of Section~\ref{nearly_holomorphic} to eliminate the non-holomorphic terms from that sequence to obtain Theorem~\ref{decomposition}.

\subsection{Frobenius Reciprocity}\label{Frobenius_reciprocity}
We review a basic construction concerning homogeneous vector bundles (cf.\@ \cite[Section~1.4.4]{CS}),
from which, in the terminology of Definition~\ref{covariant_diffop}, we obtain covariant differential operators
of order zero.

Throughout this section, let $(V, \eta)$ and~$(W, \vartheta)$ be two factors of automorphy.
Note that~$V$ and~$W$ carry the structure of a representation of $K$
by mapping $k \in K$ to $\eta(k, (i, 0))$, similarly for~$W$. We write $kv$ for $\eta(k, (i, 0))\, v$. Denote by $\Hom_K(V, W)$ the space of all $\CC$-linear maps $f \colon V \to W$ satisfying $f(kv) = kf(v)$ for all $v\in V$ and $k\in K$.

\begin{lemma}
\label{la:k_hom_to_jac_hom}
	For $f \in \Hom_K(V, W)$, the function
	\begin{gather*}
		\tilde f\colon \Jac_{1, h}(\RR)\times V \longrightarrow W, \quad (g, v) \longmapsto \vartheta(g, (i, 0))\;f\bigl(\eta(g, (i, 0))^{-1}\, v\bigr)
	\end{gather*}
	is right $K$-invariant in the first argument, i.e., $\tilde{f}(gk, v) = \tilde{f}(g, v)$ for all $k \in K$.

In particular, $\tilde{f}$ factors over the projection from $\Jac_{1, h}(\RR)\times V$ to~$\HH_{1, h}\times V$.
\end{lemma}

\begin{proof}
The calculation
\begin{align*}
  \tilde f(gk, v)
&=
  \vartheta(gk, (i, 0))\,
  f\bigl( \eta(gk, (i, 0))^{-1}\; v \bigr)
\\
&=
  \vartheta(g, (i, 0))\, k\,
  f\bigl( k^{-1}\, \eta(g, (i, 0))^{-1}\; v \bigr)
\\
&=
  \vartheta(g, (i, 0))\,
  f\bigl( \eta(g, (i, 0))^{-1}\, v \bigr)
=
  \tilde{f}(g, v)
\end{align*}
proves the claim.
\end{proof}

Given~$f \in \Hom_K(V,W)$, we continue to denote the induced map $\HH_{1, h} \times V \to W$ as in Lemma~\ref{la:k_hom_to_jac_hom} by $\tilde f$.

Conversely, given a function $F \colon \HH_{1, h}\times V \to W$, we define a map from $V$\nbd{}valued to $W$\nbd{}valued functions
\begin{gather*}
  F_* \colon \rmC^\infty(\HH_{1, h})\otimes V \longrightarrow \rmC^\infty(\HH_{1, h})\otimes W
\end{gather*}
by
\begin{gather*}
  F_*(\phi)(\tau, z) = F\bigl( (\tau, z), \phi(\tau, z) \bigr)
\end{gather*}
for $\phi \in \rmC^\infty(\HH_{1, h})\otimes V$.

To state the next theorem, recall our notation from Definition~\ref{covariant_diffop} for spaces of covariant differential operators and put
\begin{gather*}
  \Hom_{\Jac_{1, h}(\RR)}(V, W)
=
  \bigl\{
  F \in \rmC^\infty(\HH_{1, h}\times V) \otimes W
  \mathrel{\mid}
  F_* \in \Cov(\eta, \vartheta)
  \bigr\}
\tx{.}
\end{gather*}

\begin{theorem}[Frobenius Reciprocity]\label{Frobenius}
The maps
\begin{alignat*}{2}
  \Hom_K(V, W) &\longrightarrow \Hom_{\Jac_{1, h}(\RR)}(V, W)
\tx{,}\quad
&
  f &\longmapsto \tilde f
\quad\tx{and}
\\
  \Hom_{\Jac_{1, h}(\RR)}(V, W) &\longrightarrow \Hom_K(V, W)
\tx{,}\quad
&
  F &\longmapsto F\bigl( (i, 0), \;\cdot\; \bigr)
\end{alignat*}
are well-defined and mutually inverse to each other.
\end{theorem}
\begin{proof}
For $(\tau, z) = h\act{i, 0}$, we show that $\tilde{f} \in \Hom_{\Jac_{1, h}(\RR)}(V, W)$ by calculating
\begin{align*}
&
  \bigl( \tilde{f}_*(\phi\sla{\eta}{g}) \bigr) (\tau, z)
\\
={}&
  \vartheta(h, (i, 0))\,
  f\bigl( \eta(h, (i, 0))^{-1}\, \eta(g, (\tau, z))^{-1}\, \phi(g\act{\tau, z} \bigr)
\\
={}&
  \vartheta(g, (\tau, z))^{-1}\,
  \vartheta(gh, (i, 0))\,
  f\bigr( \eta(gh, (i, 0))^{-1}\, \phi(gh\act{i, 0}) \bigr)
\\
={}&
  \bigl( \tilde{f}_*(\phi) \sla{\vartheta}{g} \bigr) (\tau, z)
\tx{.}
\end{align*}
We omit the remaining routine verifications.
\end{proof}

In the next lemma, we consider the relation between $V_s$ and \mbox{$V_s \mathbin{\slash} \im(i_s^1)$},
which, by Lemma~\ref{inclusion_sym_power}, are $K$-representations via $\eta_{(k, s), m}$ and $\eta_{k, m}^{\oplus \mu(s, h)}$,
and the associated spaces of nearly holomorphic functions.

\begin{lemma}\label{section_non_holomorphic}
    Consider the map
    \begin{gather*}
        \sigma \colon   V_{s} \mathbin{\slash} \im(i_s^1)    \longrightarrow V_s,    \quad
                        Y^\nu + \im(i_s^1)                   \longmapsto     Y^\nu,
    \end{gather*}
    where we write $Y^\nu = Y_1^{\nu_1}\cdots Y_h^{\nu_h}$.
    Then
    \begin{gather*}
      \sigma \in \Hom_K\big( V_s \mathbin{\slash} \im(i_s^1),\, V_s \big)
    \tx{.}
    \end{gather*}
    Further, the map $\wtd{\sigma}_*$ induced via Theorem~\ref{Frobenius}
    satisfies $(p^1_s)_* \circ \wtd{\sigma}_* = \id$
    and restricts to a map
    \begin{gather*}
        \wtd{\sigma}_* \colon \bigl(\N_{k, m}^d\bigr)^{\oplus \mu(s, h)} \longrightarrow \N^d_{(k, s), m}
    \tx{,}
    \end{gather*}
    where~$\mu(s,h)$ is given in Lemma~\ref{inclusion_sym_power}.
\end{lemma}
\begin{proof}
For $k(\theta, \kappa) \in K$ as in \eqref{k_theta_kappa}, we find
\begin{align*}
  \sigma\bigl( k(\theta, \kappa)\, \bigl( Y^\nu + \im(i_s^1) \bigr) \bigr)
&=
  \sigma\bigl( e(k\theta)\, Y^\nu + \im(i_s^1) \bigr)
\\
&=
  e(k\theta)\, Y^\nu
=
  k(\theta, \kappa)\, \sigma\bigl( Y^\nu + \im(i_s^1) \bigr)
\end{align*}
and thus $\sigma \in \Hom_K(V_s \mathbin{\slash} \im(i_s^1), V_s)$.

For $\phi = f\, Y^\nu + \im(i_s^1)$ with $f\in \N^d_{k, m}$ and $(\tau, z) = g\act{i, 0}$ with $g \in \Jac_{1, h}(\RR)$ as in \eqref{transitivity}, we calculate
\begin{align*}
  \wtd{\sigma}_*(\phi)(\tau, z)
&=
  \eta_{(k, s), m}(g, (i, 0))\, \bigl(\eta_{k, m}(g, (i, 0))^{-1}\, f(\tau, z)\,  Y^\nu \bigr)
\\
&=
  f(\tau, z)\; \sym^s\bigl(\begin{psmatrix} 1 \slash \sqrt{y} & -\alpha \\ 0 & 1\end{psmatrix}\bigr)\, Y^\nu
\\
&=
  f(\tau, z)\, (Y_1 - \alpha_1 X)^{\nu_1}\cdots (Y_h - \alpha_h X)^{\nu_h}
\tx{.}
\end{align*}
From this we see that $(p^1_s)_* \circ \wtd{\sigma}_* = \id$.

To prove the last claim, recall that $\im(i^t_s)$ consists of polynomials divisible by $X^t$. Since in the above expression~$\alpha_j$'s appear only together with~$X$, we conclude that
\begin{gather*}
  \deg\bigl( p^t_s\circ\wtd{\sigma}_*(\phi) \bigr)
<
  \deg(f) + t \leq d + t
\tx{,}
\end{gather*}
and thus~$\wtd{\sigma}_*(\phi) \in \N^d_{(k, s), m}$.
\end{proof}

\begin{corollary}
\label{cor:section_non_holomorphic}
We have a short exact sequence
\begin{gather*}
  0 \longrightarrow \N^1_{(k+1, s-1), m} \longrightarrow \N^0_{(k, s), m} \longrightarrow \H_{k, m}^{\oplus \mu(s, h)} \longrightarrow 0
\tx{,}
\end{gather*}
which splits by a covariant differential operator of order zero.
\end{corollary}

\begin{proof}
The exact sequence is a special case of~\eqref{nearly_holomorphic_short_exact} with~$d = 0$, and the splitting differential operator was given in the preceeding lemma.
\end{proof}

\subsection{Decomposition and proof of the Main Theorem}\label{decomposition_section}

We now construct higher order differential operators to project the head of the sequence in Corollary~\ref{cor:section_non_holomorphic} onto holomorphic subspaces.
The following lemmata serve as preparation.
\begin{lemma}\label{splitting_preparation}
    Consider a diagram of vector spaces
    \begin{gather*}
        \xymatrix{
                        & 0                                         & 0                                             & 0                                                             &   \\
                        & Q_1   \ar[u] \ar@{-->}@/_1pc/[d]          & Q_2   \ar[u]                                  & Q_3   \ar[u] \ar@{-->}@/^1pc/[d]                              &   \\
            0 \ar[r]    & A     \ar[u] \ar[r] \ar@{-->}@/_1pc/[d]   & B     \ar[u] \ar[r] \ar[r]\ar@{-->}@/_1pc/[l] & C     \ar[u] \ar[r] \ar@{-->}@/_1pc/[l] \ar@{-->}@/^1pc/[d]   & 0 \\
            0 \ar[r]    & X     \ar[u] \ar[r]                       & Y     \ar[u] \ar[r]                           & Z     \ar[u] \ar[r]                                           & 0 \\
            &           0       \ar[u]                              & 0     \ar[u]                                  & 0     \ar[u]                                                  &
        }
    \end{gather*}
    with exact rows and columns, and assume that the upper row and the left and right column split.
    Then both the lower row and the middle column also split.
\end{lemma}
\begin{proof}
    A retract to $X \to Y$ is given by the composition
    \begin{gather*}
        Y \longrightarrow B \longrightarrow A \longrightarrow X.
    \end{gather*}
    Thus the lower row splits and we find a section $Z \to Y$ to $Y \to Z$.

    Since~$A \to B \to C$ splits, $\id_B$ can be decomposed as a sum of morphisms~$B \to A \to B$ and~$B \to C \to B$.
    Thus a retract to $Y \to B$ is given by the sum of the compositions
    \begin{gather*}
        B \longrightarrow C \longrightarrow Z \longrightarrow Y
    \quad\tx{and}\quad
        B \longrightarrow A \longrightarrow X \longrightarrow Y
    \tx{.}
    \end{gather*}
    This shows that the middle column splits.
\end{proof}

\begin{lemma}\label{splitting}
    Consider the commutative diagram
    \begin{gather*}
        \xymatrix{
            0 \ar[r]    & \N_{(k+1, s - 1), m}^{d + 1}                           \ar[r]%
                            & \N_{(k, s), m}^d                               \ar[r]%
                            & \bigl( \N_{k, m}^d \bigr)^{\oplus \mu(s, h)}    \ar[r]%
                                & 0                                             \\
            0 \ar[r]    & \H_{(k+1, s - 1), m}                                   \ar[r]  \ar@{^(->}[u]
                            & \H_{(k, s), m}                                 \ar[r]  \ar@{^(->}[u]
                            & \H_{k, m}^{\oplus \mu(s, h)}                   \ar[r]  \ar@{^(->}[u]
                                & 0
        }
    \end{gather*}
    assembled from the short exact sequences in~\eqref{holomorphic_short_exact} and~\eqref{nearly_holomorphic_short_exact}.
    Assume that the upper row splits and the outer vertical inclusions admit retracts, all by covariant differential operators.

    Then the inclusion in the middle column also admits a retract, and the lower row splits as well, all by covariant differential operators.
\end{lemma}
\begin{proof}
    The desired covariant differential operators can mostly be constructed as in the proof of Lemma~\ref{splitting_preparation},
    since compositions and linear combinations of covariant differential operators are again covariant differential operators.

    We only need to justify that the section $\H_{k, m}^{\oplus \mu(s, h)} \to \H_{(k, s), m}$ constructed from
    the retract $R\colon \H_{(k, s), m} \to \H_{(k + 1, s - 1), m}$ may be assumed to be a differential operator as opposed to
    a general morphism of vector spaces.

    As in Remark~\ref{non-covariant_split}, we find a section by a differential operator (not necessarily covariant)
    \begin{equation*}
        D \colon \H_{k, m}^{\oplus \mu(s, h)} \to \H_{(k, s), m}\quad \text{to}\quad (p^1_s)_* \colon \H_{(k, s), m} \to \H_{k, m}^{\oplus \mu(s, h)}.
    \end{equation*}
    We conclude that $\big(\id - (i_s^1)_*\circ R\big) \circ D$ is a covariant differential operator and a section to $(p^1_s)_*$.
\end{proof}

\begin{lemma}\label{retract_nearly_holomorphic}
    If $k - d > h/2$ and~$m$ is invertible, then the inclusion
    \begin{gather*}
        \H_{(k, s), m} \longhookrightarrow \N_{(k, s), m}^d
    \end{gather*}
    admits a retract by a covariant differential operator
\end{lemma}
\begin{proof}
    Combining the short exact sequences in~\eqref{holomorphic_short_exact} and~\eqref{nearly_holomorphic_short_exact} we obtain the commutative diagram
    \begin{gather*}
        \xymatrix{
            0 \ar[r]    & \N_{(k+1, s - 1), m}^{d + 1}                           \ar[r]
                            & \N_{(k, s), m}^d                               \ar[r]  \ar@{-->}@/_1.5pc/[l]
                            & \bigl( \N_{k, m}^d \bigr)^{\oplus \mu(s, h)}    \ar[r]  \ar@{-->}@/_1.5pc/[l]   \ar@{-->}@/^1.5pc/[d]
                                & 0                                             \\
            0 \ar[r]    & \H_{(k+1, s - 1), m}                                   \ar[r]  \ar@{^(->}[u]
                            & \H_{(k, s), m}                                 \ar[r]  \ar@{^(->}[u]
                            & \H_{k, m}^{\oplus \mu(s, h)}                   \ar[r]  \ar@{^(->}[u]
                                & 0
        }
    \end{gather*}
    with exact rows. The upper row splits by Lemma~\ref{section_non_holomorphic} and the right vertical inclusion admits a retract
    by Theorem~\ref{holomorphic_projection}.

    Observe that the diagram remains valid for~$s = 0$, if the spaces on the left are replaced by~$0$. Thus for $s = 0$, the left vertical inclusion is the zero map and the claim follows
    from Lemma~\ref{splitting}. We conclude the proof by induction on $s$, where the induction hypothesis guarantees the splitting of the left vertical map and allows us to apply Lemma~\ref{splitting} again.
\end{proof}

\begin{proof}[Proof of the Theorem~\ref{decomposition}]
    Consider the diagram
    \begin{gather*}
        \xymatrix{
            0 \ar[r]    & \N_{(k+1, s-1), m}^{1}                        \ar[r]  \ar@{-->}@/_1.5pc/[d]
                            & \N_{(k, s), m}^0                          \ar[r]
                            & \H_{k, m}^{\oplus\mu(s, h)}               \ar[r]  \ar@{-->}@/_1.5pc/[l]   & 0 \\
            0 \ar[r]    & \H_{(k+1, s-1), m}                            \ar[r]  \ar@{^(->}[u]
                            & \H_{(k, s), m}                            \ar[r]  \ar@{^(->}[u]
                            & \H_{k, m}^{\oplus\mu(s, h)}               \ar[r]  \ar@{=}[u]              & 0
        }
    \end{gather*}
    with exact rows from the proof of Lemma~\ref{retract_nearly_holomorphic} with~$d = 0$. Recall that the upper horizontal projection admits a section by Lemma~\ref{section_non_holomorphic}.
    The left vertical inclusion admits a retract by Lemma~\ref{retract_nearly_holomorphic}. Thus by Lemma~\ref{splitting},
    the lower row splits. By Lemma~\ref{diffop_and_fourier_expansion}, we obtain a decomposition
    \begin{gather*}
        \rmJ_{(k, s), m}(\Gamma) \cong \rmJ_{(k + 1, s - 1), m}(\Gamma) \oplus \rmJ_{k, m}(\Gamma)^{\oplus\mu(s, h)}
    \end{gather*}
    by a covariant differential operator.
    Now the theorem follows by induction on~$s$.
\end{proof}